\documentclass[psamsfonts]{amsart}
\usepackage{times,amsmath,amsfonts,amscd,epic,eepic}

\newtheorem*{mainproposition*}{Main Proposition}
\theoremstyle{definition}
\newtheorem*{hypotheses*}{Hypotheses for the Main Proposition}
\newtheorem*{question*}{Question}

\newtheorem{definition}{Definition}
\newtheorem*{definition*}{Definition}

\theoremstyle{remark} 
\newtheorem{remark}{Remark} 
\newtheorem*{remark*}{Remark}

\newcommand{\id}{\operatorname{id}}

\newcommand{\coder}{\operatorname{Coder}}
\newcommand{\coalg}{\operatorname{Coalg}}

\newcommand{\M}{\mathcal{M}}
\newcommand{\C}{\mathbb{C}}
\newcommand{\G}{\mathcal{G}}

\newcommand{\E}{e}
\newcommand{\ba}{\mathbf{a}}

\title{Homotopy Probability Theory I}

\author{Gabriel C. Drummond-Cole}
\address{Department of Mathematics, Northwestern University 2033
  Sheridan Road Evanston, IL 60208-2730, USA}

\thanks{This material is based in part upon work
     supported by the National Science Foundation under Award
     No. DMS-1004625.}
\email{gabriel@math.northwestern.edu}

\author{Jae-Suk Park}
\address{
  \begin{enumerate}
  \item[a.] Center for Geometry and Physics, 
Institute for Basic Science (IBS), 77 Cheongam-ro, 
Nam-gu, Pohang-si, Gyeongsangbuk-do, Korea 790-784
\item[b.] Pohang University of Science and 
Technology (POSTECH), 77 Cheongam-ro, Nam-gu, 
Pohang-si, Gyeongsangbuk-do, Korea 790-784
  \end{enumerate}
 }

\thanks{This work was supported by the IBS 
(CA1305-01). This work was supported by the 
Mid-career Researcher
  Program through NRF funded by the MEST (no. 2010-0000497)}

\email{jaesuk@postech.ac.kr}

\author{John Terilla}
\address{Department of Mathematics
  Room 211 Kissena Hall, Queens College, The City University of New
  York 65-30 Kissena Blvd Queens, NY 11367-1597 USA}
\thanks{Thanks to the Simons Center for
     Geometry and Physics for providing an excellent working environment.}

\email{jterilla@qc.cuny.edu}



\begin{document}

 \maketitle

   \begin{abstract}
     This is the first of two papers that introduce a deformation
     theoretic framework to explain and broaden a link between
     homotopy algebra and probability theory.  In this paper,
     cumulants are proved to coincide with morphisms of homotopy
     algebras.  The sequel paper outlines how the framework presented
     here can assist in the development of homotopy probability
     theory, allowing the principles of derived mathematics to
     participate in classical and noncommutative probability theory.
   \end{abstract}

\section{Introduction}
There is a surprising coincidence between homotopy morphisms in 
operadic algebra and cumulant functions in probability theory.  
In an attempt to develop this coincidence into a meaninful
relationship between homotopy
algebra and probability theory, 
this paper presents a deformation theoretic framework to explain 
the coincidence.
This framework is used in a sequel paper \cite{HPT2} to
introduce a homotopy theory of probability.

In probability theory one considers an algebra of random variables and
a linear map $\E$ from this algebra to the complex numbers $\C$.
Typically, the map $\E$ does not respect the product structure; that
is for random variables $X$ and $Y$, \[\E(XY)\neq \E(X)\E(Y).\] In
fact, the failure of $\E$ to be an algebra map measures important
correlations between random variables.  For example, the bilinear map
defined on the space of random variables by
\[\kappa_2(X,Y):=\E(XY)-\E(X)\E(Y)\]
defines the covariance.  The map $\kappa_2$ fits into an infinite
hierarchy of multilinear maps $\{\kappa_n\}$ called \emph{cumulants}.
The important notion that a set of random variables $\{X_1, \ldots,
X_n\}$ be \emph{independent} is defined using the vanishing of the
cumulants.

This paper presents a mathematical framework for studying linear maps
between algebras that do not respect the products.  The framework is a
manifestation of the following idea:
\begin{quote}
  \emph{ the failure of a map to respect structure has structure, if
    you know where to look.}
\end{quote}
Specifically, Section \ref{sec:shadow} contains a construction, which
in a simple case has the following input and output:

\bigskip
\begin{description}
\item[Input]
  \noindent
  \begin{itemize}
  \item chain complexes $C=(V,d)$ and $C'=(V',d')$
  \item degree zero bilinear maps $V\times V \to V$ and $V'\times V'
    \to V'$
  \item a chain map $f:C\to C'$
  \end{itemize}
\item[Output]
  \noindent
  \begin{itemize}
  \item a sequence of degree zero multlinear maps $\{f_n:V^{\times n} \to
    V'\}_{n=1}^\infty$.
  \end{itemize}
\end{description}

For the input, the differentials and the chain map are not assumed to
have any compatibility with the bilinear maps.  The sequence of maps
$\{f_n\}$ in the output constitute an $A_\infty$ algebra morphism
between two $A_\infty$ algebras which arise during the construction.

A probability space provides an example of the input data for the
construction.  The chain complexes $C$ and $C'$ are, respectively,
$(V,0)$ and $(\C,0)$, the space of random variables and the complex
numbers both with zero differential. The bilinear maps are given by
the products of random variables and complex numbers.  The chain map
$C\to C'$ is the expectation $\E:V\to \C$.  Thus, the construction can
be applied producing an $A_\infty$ morphism
$$\{e_n:V^n\to \C\}_{n=1}^\infty$$
as output.  The main result of this paper is that this $A_\infty$ morphism
coincides with the Boolean cumulants (see Definition \ref{def:cumulants} and Remark \ref{rk:cumulants}): $\kappa_n=e_n$ for all $n$.

The paper proceeds as follows.  Section \ref{sec:shadow} contains the
details of the construction of the output $A_\infty$ morphism.  While
$A_\infty$ structures have been studied in topology since the middle
of the last century \cite{St1}, Section \ref{sec:shadow} also includes a brief overview of $A_\infty$ algebras and
morphisms.  This may be helpful to readers unfamiliar with $A_\infty$
structures, and it also highlights the notion of 
equivalence in their moduli spaces, which is an important tool in our approach.
Section \ref{sec:cumulants} contains the definition of the cumulants
and the main proposition.  
There is nothing novel in the proof of the main proposition, 
which is a formal verification.  The novelty 
is the statement, which connects the two previously unlinked concepts
of homotopy morphisms and cumulants.

The authors would like to thank Tyler Bryson, Joseph Hirsh, Tom
LeGatta, and Bruno Vallette for many helpful discussions.

\section{$A_\infty$ algebras and morphisms}\label{sec:shadow}
The book \cite{LV} is good reference for $A_\infty$ algebras.
\subsection{Definitions}
Let $V$ be a graded vector space.  Let $V^{\otimes n}$ denote the
$n$th tensor power of $V$ and $TV=\oplus_{n=1}^\infty V^{\otimes n}$.
As a direct sum, linear maps from $TV$ to a vector space $W$
correspond to collections of linear maps $\{V^{\otimes n} \to
W\}_{n=1}^\infty$.  Also, $TV$ is a coalgebra, free in a certain
sense, so that coalgebra maps from a coalgebra $\mathcal{C}$ to $TV$
correspond to linear maps $\mathcal{C}\to V$.  This freeness also
implies that coderivations from a coalgebra $\mathcal{C}$ to $TV$
correspond to linear maps $\mathcal{C}\to V$.  All these
correspondences are one-to-one: any linear map $C\to V$ can be lifted
uniquely to a coderivation $\mathcal{C} \to TV$, or lifted uniquely to
a coalgebra map $\mathcal{C} \to TV$.
\begin{definition}
  \emph{An $A_\infty$ algebra} is a pair $(V,D)$ where $V$ is a graded
  vector space and $D:TV \to TV$ is a degree one\footnote{Readers
    familiar with $A_\infty$ algebras will be aware that the
    definition of an $A_\infty$ algebra usually involves a shift of
    degree, but no degree shift is used in this paper.}
  coderivation satisfying $D^2=0$.  \emph{An $A_\infty$ morphism}
  between two $A_\infty$ algebras $(V,D)$ and $(V',D')$ is a
  differential coalgebra map $F:(TV,D)\to (TV',D')$.  In other words,
  an $A_\infty$ map from $(V,D)$ to $(V',D')$ is a degree zero
  coalgebra map $F:TV \to TV'$ satisfying $FD=D'F$.
\end{definition}

The identification
\begin{equation}
  \label{eq:coder}
  \coder(TV,TV)\simeq
  \prod_{n=1}^\infty \hom(V^{\otimes n},V)
\end{equation}
can be used to give the data of an $A_\infty$ algebra.  That is, a
coderivation $D:TV \to TV$ can be given by a sequence $\{d_n\}$ of
degree one linear maps $d_n:V^{\otimes n} \to V$ on the graded vector
space $V$.  The condition that $D^2=0$ implies an infinite number of
relations satisfied by various compositions among the $\{d_n\}$.
Likewise, the identification
\begin{equation}
  \label{eq:coalg}
  \coalg(TV,TV')\simeq \prod_{n=1}^\infty \hom(V^{\otimes n},V')
\end{equation}
can be used to give the data of an $A_\infty$ morphism between $(V,D)$
and $(V',D')$.  That is, a coalgebra map $F:TV \to TV'$ can be given
by a sequence $\{f_n\}$ of degree zero maps $f_n:V^{\otimes} \to V'$.
The condition that $FD=D' F$ encodes an infinite number of relations
among the $\{f_n\}$, the $\{d_n\}$ and the $\{d'_n\}$.  So, the
identifications in Equations \eqref{eq:coder} and \eqref{eq:coalg}
provide two equivalent ways of describing $A_\infty$ algebras and
morphisms and it is convenient to move between the two ways since certain
notions or computations are easier to describe in one or the the other
description of the equivalent data.  For example, it is
straightforward to compose two differential coalgebra maps
$F:(TV,D)\to (TV',D')$ and $G:(TV',D')\to (TV'',D'')$ as $GF:(TV,D)\to
(TV'',D'')$ and thus define the composition of $A_\infty$ morphisms,
but it is more involved to express the components $(gf)_n:V^{\otimes
  n} \to V''$ in terms of the $f_k:V^{\otimes k} \to V'$ and
$g_m:V'^{\otimes m} \to V''$.

Note that $A_\infty$ structures can be transported via isomorphisms.
In particular, if $(V,D)$ is an $A_\infty$ algebra and $G:TV \to TW$
is any degree zero coalgebra isomorphism, then for \[D^G:=G^{-1}DG,\]
the pair $(W,D^G)$ is again an $A_\infty$ algebra.

Morphisms can be transported as well.  If $F$ is an $A_\infty$
morphism between $(V,D)$ and $(V', D')$ and $G:TV \to TW$ and $H:TV'
\to TW'$ are coalgebra isomorphisms, then \[F^{G,H}:=H^{-1}FG\]
is an $A_\infty$ morphism between the $A_\infty$ algebras with the
transported structures $(W,D^G)$
and $(W',(D')^H)$.

\begin{definition}
  Two $A_\infty$ algebras $(V,D)$ and $(V',D')$ are \emph{equivalent}
  if there exists a coalgebra isomorphism $G:TV \to TV'$ so that
  $D'=D^G.$
\end{definition}

\subsection{Spaces of $A_\infty$ algebras}

Consider an $A_\infty$ algebra $(V,D)$.  The condition that $D$ has
degree one and that $D^2=0$ imply that the first component $d_1:V \to
V$ of $D$ has degree one and satisfies $d_1\circ d_1=0$.  So, the pair
$(V,d_1)$ is a chain complex.  Often it is appropriate to view the
chain complex $(V,d_1)$ as a fundamental object and to consider the
remaining components $d_2, d_3, \ldots$ of $D$ as structure on the
chain complex $(V,d_1)$.

\begin{definition}
  Let $C=(V,d)$ be a chain complex.  An \emph{$A_\infty$ structure on
    $C$} is an $A_\infty$ algebra $(V,D)$ with $d_1=d.$ Let $\M_C$
  denote the set of $A_\infty$ structures on $C$.
\end{definition}

\begin{definition}\label{gauge_group_def}
  Let $C=(V,d)$ be a chain complex.  The \emph{gauge group} $\G_C$ is
  the subgroup of $\mathrm{GL}(TV)$ consisting of degree zero
  coalgebra automorphisms with first component $\id:V\to V$.  The
  gauge group $\G_C$ acts, on the right, by conjugation on
  $\M_C$: \[D\cdot G = D^G.\]
\end{definition}

\subsection{The gauge group $\G_C$}\label{sec:gauge}

The \emph{gauge group} $\G_C$ as defined in Definition
\ref{gauge_group_def} is a Lie subgroup of $\mathrm{GL}(TV)$.
The Lie algebra of $\G_C$ is the Lie subalgebra of $\mathrm{gl}(TV)$
consisting of all degree zero coderivations $TV \to TV$ with first
component $0:V\to V$. This Lie subalgebra can be identified via
Equation \eqref{eq:coder} with the vector space of degree zero maps
$\prod_{k=2}^\infty \hom(V^{\otimes k},V)$.  Any map $a:V^{\otimes
  k}\to V$ with $k>1$ can be lifted to a coderivation $A:TV\to TV$ and
then exponentiated to obtain a gauge group element $\ba \in\M_C$ where
$$\ba=\exp(A):=\id+A+\frac{A^2}{2!}+\frac{A^3}{3!}+\cdots$$
The orbits of the one parameter subgroup $t\mapsto \exp(t A)$ of the
gauge group $\G_C$ are curves in the space $\M_C$.  Each of these
curves connects an $A_\infty$ structure $D$ at $t=0$ to an equivalent
$A_\infty$ structure $D^{\ba}$ at $t=1$.
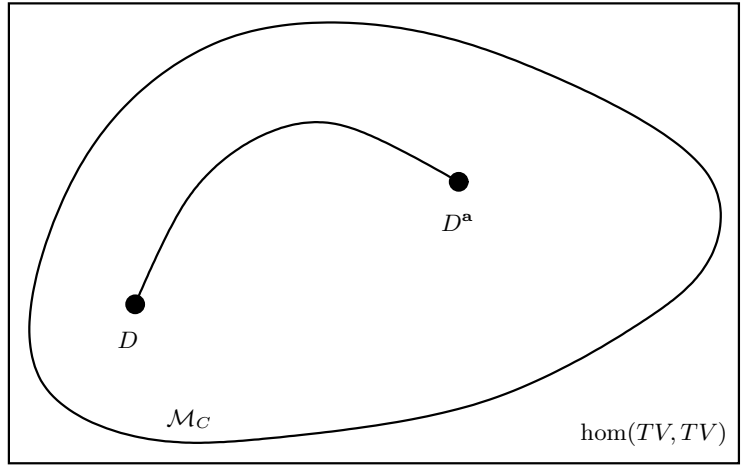
\begin{figure}[h!]
\bigskip
 \begin{center}
\setlength{\unitlength}{0.0003in}
{\renewcommand{\dashlinestretch}{30}
\begin{picture}(11926,7297)(0,-10)
\thicklines
\path(7800,746)(7945,793)(8087,843)
	(8226,895)(8362,948)(8495,1003)
	(8623,1058)(8748,1113)(8869,1169)
	(8986,1225)(9099,1281)(9208,1336)
	(9313,1390)(9414,1444)(9512,1496)
	(9605,1548)(9695,1599)(9781,1648)
	(9864,1696)(9944,1743)(10020,1789)
	(10093,1833)(10164,1876)(10231,1918)
	(10296,1959)(10359,1999)(10419,2037)
	(10477,2074)(10533,2111)(10587,2146)
	(10640,2181)(10690,2215)(10740,2248)
	(10787,2281)(10834,2313)(10879,2345)
	(10923,2377)(10966,2408)(11009,2440)
	(11050,2471)(11091,2502)(11131,2534)
	(11170,2566)(11209,2599)(11247,2632)
	(11285,2666)(11322,2700)(11358,2736)
	(11395,2772)(11431,2810)(11466,2849)
	(11501,2889)(11535,2930)(11568,2973)
	(11601,3018)(11633,3064)(11665,3113)
	(11695,3163)(11724,3215)(11752,3269)
	(11779,3325)(11803,3383)(11826,3443)
	(11847,3506)(11866,3571)(11882,3637)
	(11895,3706)(11905,3777)(11912,3850)
	(11914,3925)(11913,4001)(11907,4079)
	(11896,4158)(11879,4238)(11857,4319)
	(11828,4401)(11792,4483)(11750,4565)
	(11700,4646)(11650,4716)(11596,4784)
	(11536,4852)(11471,4919)(11403,4985)
	(11332,5049)(11257,5112)(11180,5174)
	(11100,5234)(11019,5294)(10936,5351)
	(10851,5408)(10766,5463)(10681,5517)
	(10594,5569)(10508,5620)(10422,5669)
	(10335,5718)(10249,5765)(10164,5811)
	(10079,5855)(9994,5899)(9911,5941)
	(9828,5983)(9746,6023)(9665,6062)
	(9584,6100)(9505,6138)(9426,6174)
	(9348,6210)(9271,6245)(9195,6279)
	(9120,6313)(9045,6345)(8971,6378)
	(8897,6409)(8825,6440)(8752,6471)
	(8680,6501)(8608,6530)(8536,6559)
	(8465,6588)(8393,6616)(8322,6644)
	(8250,6671)(8178,6698)(8106,6725)
	(8033,6751)(7959,6777)(7885,6803)
	(7810,6828)(7735,6853)(7658,6878)
	(7580,6902)(7501,6926)(7421,6950)
	(7339,6973)(7256,6995)(7171,7018)
	(7084,7039)(6996,7060)(6906,7081)
	(6813,7101)(6719,7120)(6623,7139)
	(6524,7156)(6423,7173)(6320,7189)
	(6215,7203)(6107,7217)(5998,7229)
	(5885,7240)(5771,7249)(5654,7257)
	(5535,7263)(5414,7267)(5290,7270)
	(5165,7270)(5038,7267)(4910,7263)
	(4780,7255)(4649,7245)(4518,7232)
	(4385,7216)(4253,7197)(4120,7174)
	(3988,7148)(3857,7118)(3728,7084)
	(3600,7046)(3469,7002)(3340,6954)
	(3215,6902)(3093,6847)(2974,6788)
	(2859,6727)(2747,6664)(2639,6599)
	(2535,6533)(2434,6464)(2337,6395)
	(2244,6325)(2154,6254)(2068,6183)
	(1985,6112)(1905,6040)(1829,5969)
	(1756,5897)(1687,5826)(1620,5756)
	(1556,5685)(1495,5615)(1436,5546)
	(1381,5477)(1327,5409)(1276,5342)
	(1227,5275)(1180,5209)(1135,5143)
	(1091,5078)(1050,5014)(1010,4950)
	(972,4887)(934,4824)(899,4761)
	(864,4699)(831,4638)(799,4576)
	(767,4515)(737,4454)(707,4393)
	(678,4332)(650,4271)(622,4210)
	(595,4149)(569,4087)(543,4025)
	(517,3963)(492,3900)(467,3837)
	(443,3773)(419,3709)(395,3644)
	(372,3578)(349,3511)(326,3443)
	(304,3375)(282,3305)(261,3234)
	(240,3163)(220,3090)(200,3016)
	(181,2941)(163,2865)(146,2787)
	(130,2709)(115,2629)(101,2548)
	(89,2466)(78,2383)(69,2300)
	(62,2215)(57,2130)(54,2044)
	(54,1958)(57,1871)(63,1785)
	(73,1698)(85,1612)(102,1527)
	(123,1443)(148,1359)(179,1278)
	(214,1198)(254,1121)(300,1046)
	(350,976)(406,909)(466,845)
	(530,784)(597,726)(668,671)
	(741,618)(817,569)(894,523)
	(973,480)(1053,439)(1134,401)
	(1215,366)(1297,333)(1379,302)
	(1461,274)(1542,248)(1623,224)
	(1703,202)(1783,182)(1862,164)
	(1940,147)(2017,133)(2094,119)
	(2169,108)(2243,97)(2316,88)
	(2388,80)(2460,74)(2530,68)
	(2599,64)(2668,60)(2736,57)
	(2803,55)(2869,54)(2935,54)
	(3000,54)(3065,55)(3129,56)
	(3193,58)(3258,61)(3322,64)
	(3386,67)(3450,71)(3515,75)
	(3579,80)(3645,85)(3711,90)
	(3778,95)(3845,101)(3914,107)
	(3984,113)(4055,120)(4127,127)
	(4200,134)(4275,141)(4352,149)
	(4431,157)(4511,166)(4594,174)
	(4678,184)(4765,193)(4854,203)
	(4945,213)(5039,224)(5135,236)
	(5234,247)(5335,260)(5439,273)
	(5545,287)(5655,302)(5766,318)
	(5881,334)(5997,352)(6117,370)
	(6238,390)(6362,411)(6488,433)
	(6615,457)(6745,482)(6875,509)
	(7007,537)(7139,567)(7272,599)
	(7405,633)(7538,668)(7670,706)(7800,746)
\thicklines
\path(1866,2432)(1867,2433)(1868,2436)
	(1870,2442)(1874,2451)(1880,2464)
	(1887,2481)(1896,2503)(1908,2530)
	(1922,2562)(1939,2600)(1957,2642)
	(1978,2689)(2001,2741)(2025,2797)
	(2052,2856)(2080,2918)(2109,2982)
	(2139,3048)(2170,3116)(2201,3184)
	(2233,3252)(2265,3320)(2297,3387)
	(2329,3453)(2360,3518)(2392,3581)
	(2423,3642)(2453,3702)(2483,3760)
	(2513,3816)(2542,3870)(2570,3922)
	(2599,3972)(2627,4020)(2654,4067)
	(2682,4112)(2709,4155)(2737,4197)
	(2764,4238)(2791,4278)(2819,4316)
	(2847,4354)(2875,4391)(2903,4427)
	(2932,4462)(2961,4497)(2991,4532)
	(3023,4568)(3055,4603)(3088,4638)
	(3122,4673)(3156,4708)(3192,4742)
	(3228,4776)(3265,4810)(3303,4844)
	(3342,4877)(3382,4910)(3423,4943)
	(3465,4976)(3508,5008)(3551,5039)
	(3596,5071)(3641,5101)(3687,5131)
	(3733,5161)(3781,5189)(3828,5217)
	(3877,5244)(3926,5270)(3975,5296)
	(4024,5320)(4074,5343)(4124,5365)
	(4174,5386)(4224,5406)(4274,5425)
	(4324,5442)(4373,5459)(4423,5474)
	(4472,5488)(4521,5500)(4570,5512)
	(4618,5522)(4666,5530)(4714,5538)
	(4761,5544)(4808,5549)(4855,5553)
	(4902,5556)(4948,5558)(4995,5558)
	(5041,5557)(5084,5555)(5126,5552)
	(5169,5548)(5213,5543)(5256,5536)
	(5301,5529)(5345,5520)(5391,5509)
	(5437,5497)(5485,5484)(5533,5470)
	(5583,5454)(5634,5436)(5687,5417)
	(5741,5396)(5797,5374)(5854,5350)
	(5913,5325)(5974,5297)(6037,5269)
	(6101,5239)(6168,5207)(6235,5174)
	(6304,5139)(6374,5104)(6445,5067)
	(6517,5030)(6589,4992)(6660,4954)
	(6731,4916)(6801,4878)(6869,4841)
	(6935,4805)(6998,4770)(7057,4737)
	(7113,4705)(7165,4676)(7211,4650)
	(7253,4626)(7290,4605)(7321,4587)
	(7348,4572)(7369,4559)(7386,4550)
	(7398,4543)(7416,4532)
\blacken\path(7297.963,4568.976)(7416.000,4532.000)(7329.250,4620.172)(7297.963,4568.976)
\put(1566,1682){\makebox(0,0)[lb]{\small $D$}}
\put(1866,2432){\blacken\ellipse{300}{300}}
\put(7116,3682){\makebox(0,0)[lb]{\small $D^{\ba}$}}
\put(7416,4532){\blacken\ellipse{300}{300}}

%
\path(-300,-300)(12226,-300)(12226,7597)(-300,7597)(-300,-300)
\put(12026,0){\makebox(0,0)[rb]{\small $\hom(TV,TV)$}}
\put(2400,300){\makebox(0,0)[lb]{\small $\M_C$}}

\end{picture}
}
 \end{center}
\bigskip
 \caption{A picture of the family of $A_\infty$ structures $\M_C$ on
    a chain complex $C=(V,d)$ as a subset of $\hom(TV,TV)$.
   There is a basepoint $D:TV\to TV$ which is the lift of $d:V\to V$
   as a coderivation---the basepoint corresponds to $C$ itself with
   zero additional structure.  
   The path is part of the orbit of the one parameter subgroup of the
   gauge group passing through the basepoint $D$ at $t=0$ and an
   equivalent $A_\infty$ structure $D^{\ba}:TV \to TV$ at $t=1$. }

\end{figure}

\begin{remark} In classical deformation theory, equivalent structures
  are identified to form a quotient moduli set.  Rather than
  identifying equivalent structures, a simplicial moduli space can be
  constructed.  The points of the simplicial moduli space consist of
  all structures.  The paths in the simplicial moduli space consist of
  equivalences between structures.  Higher dimensional parts
  correspond to equivalences between equivalences.  This paper
  involves $A_\infty$ structures which are equivalent to trivial
  $A_\infty$ structures.  In order to see the application to
  probability theory, gauge equivalent structures should not be
  identified, so the relevant moduli space is the simplicial moduli
  space.
\end{remark}

\subsection{Spaces of $A_\infty$ morphisms}
Consider an $A_\infty$ morphism $F:TV\to TV'$ between $A_\infty$
algebras $(V,D)$ and $(V',D')$.  The conditions that $F$ has degree
zero and that $FD=D' F$ imply that the first component $f_1:V\to V'$
has degree zero and satisfies $f_1d_1=d'_1 f_1$.  Thus $f_1:(V,d_1)\to
(V',d'_1)$ is a chain map.  Here, the chain map $f_1:(V,d_1)\to
(V',d'_1)$ is viewed as a fundamental object and the remaining
components $f_2, f_3, \ldots$ of $F$ are viewed as a structure on the
chain map $f_1$.

\begin{definition}
  Let $C=(V,d)$ and $C'=(V',d')$ be chain complexes and let $f:C\to
  C'$ be a chain map.  An \emph{$A_\infty$ morphism on $f$} is an $A_\infty$
  morphism $F:(TV,D)\to (TV,D')$ between two $A_\infty$ structures on
  $C$ and $C'$ with $f_1=f.$  Let $$\M_f\subset\hom(TV,TV)\times
  \hom(TV,TV')\times \hom(TV',TV')$$ be the set of $A_\infty$
  morphisms on $f$ between $A_\infty$ structures on $C$ and $C'$.
  That is, a triple $(D,F,D')\in \M_f$ consists of degree one, square
  zero coderivations $D:TV\to TV$ and $D':TV'\to TV'$ with $d_1=d$ and
  $d'_1=d'$ and a degree zero coalgebra map $F:TV \to TV'$ satisfying
  $FD=D'F$ with $f_1=f$.

  The product of gauge groups $\G_C\times \G_{C'}$ acts, on the right,
  on $\M_f$ by $$ (D,F,D')\cdot (G,H):=(D^G,F^{H,G},(D')^H).$$
\end{definition}

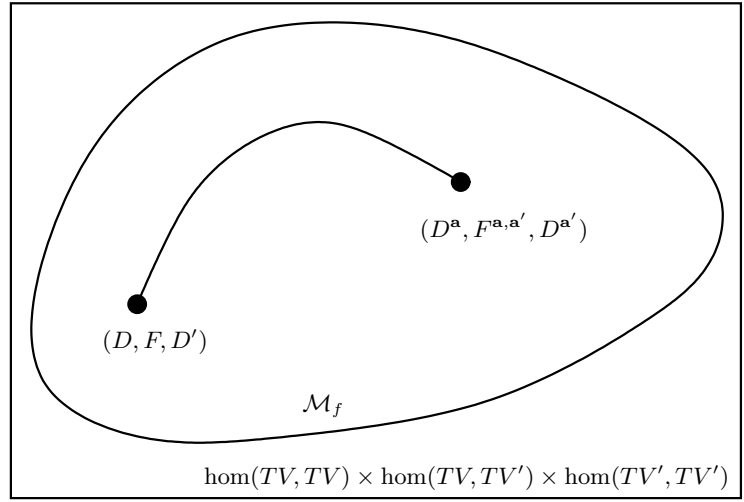
\begin{figure}
\bigskip
  \begin{center}
\setlength{\unitlength}{0.0003in}
{\renewcommand{\dashlinestretch}{30}
\begin{picture}(11926,7297)(0,-10)
\thicklines
\path(7800,746)(7945,793)(8087,843)
	(8226,895)(8362,948)(8495,1003)
	(8623,1058)(8748,1113)(8869,1169)
	(8986,1225)(9099,1281)(9208,1336)
	(9313,1390)(9414,1444)(9512,1496)
	(9605,1548)(9695,1599)(9781,1648)
	(9864,1696)(9944,1743)(10020,1789)
	(10093,1833)(10164,1876)(10231,1918)
	(10296,1959)(10359,1999)(10419,2037)
	(10477,2074)(10533,2111)(10587,2146)
	(10640,2181)(10690,2215)(10740,2248)
	(10787,2281)(10834,2313)(10879,2345)
	(10923,2377)(10966,2408)(11009,2440)
	(11050,2471)(11091,2502)(11131,2534)
	(11170,2566)(11209,2599)(11247,2632)
	(11285,2666)(11322,2700)(11358,2736)
	(11395,2772)(11431,2810)(11466,2849)
	(11501,2889)(11535,2930)(11568,2973)
	(11601,3018)(11633,3064)(11665,3113)
	(11695,3163)(11724,3215)(11752,3269)
	(11779,3325)(11803,3383)(11826,3443)
	(11847,3506)(11866,3571)(11882,3637)
	(11895,3706)(11905,3777)(11912,3850)
	(11914,3925)(11913,4001)(11907,4079)
	(11896,4158)(11879,4238)(11857,4319)
	(11828,4401)(11792,4483)(11750,4565)
	(11700,4646)(11650,4716)(11596,4784)
	(11536,4852)(11471,4919)(11403,4985)
	(11332,5049)(11257,5112)(11180,5174)
	(11100,5234)(11019,5294)(10936,5351)
	(10851,5408)(10766,5463)(10681,5517)
	(10594,5569)(10508,5620)(10422,5669)
	(10335,5718)(10249,5765)(10164,5811)
	(10079,5855)(9994,5899)(9911,5941)
	(9828,5983)(9746,6023)(9665,6062)
	(9584,6100)(9505,6138)(9426,6174)
	(9348,6210)(9271,6245)(9195,6279)
	(9120,6313)(9045,6345)(8971,6378)
	(8897,6409)(8825,6440)(8752,6471)
	(8680,6501)(8608,6530)(8536,6559)
	(8465,6588)(8393,6616)(8322,6644)
	(8250,6671)(8178,6698)(8106,6725)
	(8033,6751)(7959,6777)(7885,6803)
	(7810,6828)(7735,6853)(7658,6878)
	(7580,6902)(7501,6926)(7421,6950)
	(7339,6973)(7256,6995)(7171,7018)
	(7084,7039)(6996,7060)(6906,7081)
	(6813,7101)(6719,7120)(6623,7139)
	(6524,7156)(6423,7173)(6320,7189)
	(6215,7203)(6107,7217)(5998,7229)
	(5885,7240)(5771,7249)(5654,7257)
	(5535,7263)(5414,7267)(5290,7270)
	(5165,7270)(5038,7267)(4910,7263)
	(4780,7255)(4649,7245)(4518,7232)
	(4385,7216)(4253,7197)(4120,7174)
	(3988,7148)(3857,7118)(3728,7084)
	(3600,7046)(3469,7002)(3340,6954)
	(3215,6902)(3093,6847)(2974,6788)
	(2859,6727)(2747,6664)(2639,6599)
	(2535,6533)(2434,6464)(2337,6395)
	(2244,6325)(2154,6254)(2068,6183)
	(1985,6112)(1905,6040)(1829,5969)
	(1756,5897)(1687,5826)(1620,5756)
	(1556,5685)(1495,5615)(1436,5546)
	(1381,5477)(1327,5409)(1276,5342)
	(1227,5275)(1180,5209)(1135,5143)
	(1091,5078)(1050,5014)(1010,4950)
	(972,4887)(934,4824)(899,4761)
	(864,4699)(831,4638)(799,4576)
	(767,4515)(737,4454)(707,4393)
	(678,4332)(650,4271)(622,4210)
	(595,4149)(569,4087)(543,4025)
	(517,3963)(492,3900)(467,3837)
	(443,3773)(419,3709)(395,3644)
	(372,3578)(349,3511)(326,3443)
	(304,3375)(282,3305)(261,3234)
	(240,3163)(220,3090)(200,3016)
	(181,2941)(163,2865)(146,2787)
	(130,2709)(115,2629)(101,2548)
	(89,2466)(78,2383)(69,2300)
	(62,2215)(57,2130)(54,2044)
	(54,1958)(57,1871)(63,1785)
	(73,1698)(85,1612)(102,1527)
	(123,1443)(148,1359)(179,1278)
	(214,1198)(254,1121)(300,1046)
	(350,976)(406,909)(466,845)
	(530,784)(597,726)(668,671)
	(741,618)(817,569)(894,523)
	(973,480)(1053,439)(1134,401)
	(1215,366)(1297,333)(1379,302)
	(1461,274)(1542,248)(1623,224)
	(1703,202)(1783,182)(1862,164)
	(1940,147)(2017,133)(2094,119)
	(2169,108)(2243,97)(2316,88)
	(2388,80)(2460,74)(2530,68)
	(2599,64)(2668,60)(2736,57)
	(2803,55)(2869,54)(2935,54)
	(3000,54)(3065,55)(3129,56)
	(3193,58)(3258,61)(3322,64)
	(3386,67)(3450,71)(3515,75)
	(3579,80)(3645,85)(3711,90)
	(3778,95)(3845,101)(3914,107)
	(3984,113)(4055,120)(4127,127)
	(4200,134)(4275,141)(4352,149)
	(4431,157)(4511,166)(4594,174)
	(4678,184)(4765,193)(4854,203)
	(4945,213)(5039,224)(5135,236)
	(5234,247)(5335,260)(5439,273)
	(5545,287)(5655,302)(5766,318)
	(5881,334)(5997,352)(6117,370)
	(6238,390)(6362,411)(6488,433)
	(6615,457)(6745,482)(6875,509)
	(7007,537)(7139,567)(7272,599)
	(7405,633)(7538,668)(7670,706)(7800,746)

\thicklines

\path(1866,2432)(1867,2433)(1868,2436)
	(1870,2442)(1874,2451)(1880,2464)
	(1887,2481)(1896,2503)(1908,2530)
	(1922,2562)(1939,2600)(1957,2642)
	(1978,2689)(2001,2741)(2025,2797)
	(2052,2856)(2080,2918)(2109,2982)
	(2139,3048)(2170,3116)(2201,3184)
	(2233,3252)(2265,3320)(2297,3387)
	(2329,3453)(2360,3518)(2392,3581)
	(2423,3642)(2453,3702)(2483,3760)
	(2513,3816)(2542,3870)(2570,3922)
	(2599,3972)(2627,4020)(2654,4067)
	(2682,4112)(2709,4155)(2737,4197)
	(2764,4238)(2791,4278)(2819,4316)
	(2847,4354)(2875,4391)(2903,4427)
	(2932,4462)(2961,4497)(2991,4532)
	(3023,4568)(3055,4603)(3088,4638)
	(3122,4673)(3156,4708)(3192,4742)
	(3228,4776)(3265,4810)(3303,4844)
	(3342,4877)(3382,4910)(3423,4943)
	(3465,4976)(3508,5008)(3551,5039)
	(3596,5071)(3641,5101)(3687,5131)
	(3733,5161)(3781,5189)(3828,5217)
	(3877,5244)(3926,5270)(3975,5296)
	(4024,5320)(4074,5343)(4124,5365)
	(4174,5386)(4224,5406)(4274,5425)
	(4324,5442)(4373,5459)(4423,5474)
	(4472,5488)(4521,5500)(4570,5512)
	(4618,5522)(4666,5530)(4714,5538)
	(4761,5544)(4808,5549)(4855,5553)
	(4902,5556)(4948,5558)(4995,5558)
	(5041,5557)(5084,5555)(5126,5552)
	(5169,5548)(5213,5543)(5256,5536)
	(5301,5529)(5345,5520)(5391,5509)
	(5437,5497)(5485,5484)(5533,5470)
	(5583,5454)(5634,5436)(5687,5417)
	(5741,5396)(5797,5374)(5854,5350)
	(5913,5325)(5974,5297)(6037,5269)
	(6101,5239)(6168,5207)(6235,5174)
	(6304,5139)(6374,5104)(6445,5067)
	(6517,5030)(6589,4992)(6660,4954)
	(6731,4916)(6801,4878)(6869,4841)
	(6935,4805)(6998,4770)(7057,4737)
	(7113,4705)(7165,4676)(7211,4650)
	(7253,4626)(7290,4605)(7321,4587)
	(7348,4572)(7369,4559)(7386,4550)
	(7398,4543)(7416,4532)
\blacken\path(7297.963,4568.976)(7416.000,4532.000)(7329.250,4620.172)(7297.963,4568.976)
\put(1266,1582){\makebox(0,0)[lb]{\small $(D,F,D')$}}
\put(1866,2432){\blacken\ellipse{300}{300}}
\put(6716,3482){\makebox(0,0)[lb]{\small
    $(D^\ba,F^{\ba,\ba'},D^{\ba'})$}}
\put(7416,4532){\blacken\ellipse{300}{300}}

\path(-300,-900)(12226,-900)(12226,7597)(-300,7597)(-300,-900)
\put(12026,-750){\makebox(0,0)[rb]{\small $\hom(TV,TV)\times \hom(TV,TV')\times \hom(TV',TV')$}}
\put(4700,500){\makebox(0,0)[lb]{\small $\M_f$}}

\end{picture}
}
  \end{center}
\bigskip
\caption{The picture for morphisms equivalent to chain maps is
    similar to that of structures equivalent to chain complexes.  }
\end{figure}

\begin{remark}\label{rm:trivial}
  A seemingly trivial situation will be important in the next section.
  An ungraded vector space $V$ can be considered a chain complex
  $C=(V,0)$ by setting the degree of every element of $V$ to be zero
  and setting the differential $d=0$.  The gauge group $\G_C$ acts
  trivially on $\M_C$ when $d=0$ since $G^{-1} 0 G = 0$ for any
  $G\in\M_C$.  Any linear map $f:V\to V'$ between vector spaces $V$
  and $V'$ is a chain map between $C=(V,0)$ and $C'=(V',0)$ when $V$
  and $V'$ are regarded as chain complexes with zero differentials.
  If $a:V^{\otimes k} \to V$ and $a':V'^{\otimes m} \to V'$, then $a,
  a'$ can be lifted to coderivations $A, A'$ and exponentiated to
  obtain coalgebra automorphisms $\ba \in \G_C$ and $\ba' \in
  \G_{C'}$.  The $A_\infty$ structures $D^{\ba}$ and $D'^{\ba '}$ on
  $C$ and $C'$ are identically zero, but the morphism $F^{\ba,\ba'}:
  TV\to TV'$ is typically nonzero.  That is, $F^{\ba,\ba'}$ is a
  nonzero $A_\infty$ morphism between two zero $A_\infty$ structures.

  For example, a straightforward computation in the simple case that
  $a:V^{\otimes 2} \to V$ and $a':V'^{\otimes 2} \to V$ reveals the
  first two structure morphisms $f_1^{\ba,\ba'}:V\to V'$ and
  $f_2^{\ba,\ba'}:V^{\otimes 2} \to V'$ to be the maps defined by
  \begin{gather*}
    f_1^{\ba,\ba'}(X)=f(X)\\
    f_2^{\ba,\ba'}(X_1,X_2)=f(a(X_1, X_2))-a'(f(X_1),f(X_2)).
  \end{gather*}
\end{remark}

\section{Probability spaces and cumulants}\label{sec:cumulants}
\subsection{Probability spaces}
One modern approach to probability theory (see, for example,
\cite{TAO2}) begins with the following definition:
\begin{definition}
  A \emph{probability space} is a triple $(V,\E,a)$ where $V$ is a
  complex vector space, $\E:V\to \C$ is a linear function, and
  $a:V\times V \to V$ is an associative bilinear product.  Elements of
  $V$ are called \emph{random variables} and the number $\E(X)$ is
  called \emph{the expected value} of the random variable $X\in V$.
  The notation $X_1X_2$ is used for $a(X_1,X_2)$. 
  Multiplication of complex numbers is denoted by $a'$. 
  For $n>1$, the notation $\alpha_n$ (and
  $\alpha'_n$) is used for the linear maps $V^{\otimes n}\to V$
  (and $\C^{\otimes n}\to \C$) obtained by repeated multiplication;
  for $n=1$, $\alpha_1=\id_V$ (and $\alpha'_1=\id_\C$).
  The expectation values of products $\E(X_1\cdots X_n)$ are called
  \emph{joint moments}.
\end{definition}
The product on $V$ is not assumed to be commutative.  Elements of $V$
are sometimes called \emph{observables} when the product is not
commutative, but here no special terminology is used to distinguish
between commutative and noncommutative probability spaces.

\subsection{Cumulants}
\begin{definition}\label{def:cumulants}
  Let $(V,\E,a)$ be a probability space.  The \emph{$n$th cumulant} of
  $(V,\E,a)$ is the linear map $\kappa_n:V^{\otimes n} \to \C$ defined
  recursively by the following equation:
  \begin{equation}\label{eq:cumulants}
    \E\circ \alpha_n =
    \sum_{k=1}^n \alpha'_k\circ\left( \sum_P \bigotimes_{j=1}^k \kappa_{n_j}\right)
  \end{equation}
  where $P$ ranges over all ordered partitions $P=(n_1,\ldots, n_k)$
  with $\displaystyle\sum_{j=1}^k n_j = n$.
\end{definition}
\begin{remark}\label{rk:cumulants}
  The cumulants defined above have been called \emph{Boolean
    cumulants} \cite{SW} to distinguish these
  cumulants from the \emph{classical cumulants}, which are defined in
  the special case when the product on $V$ is commutative, and the
  \emph{free cumulants} which are important in free probability theory
  \cite{NSp}.  For a survey of various kinds of cumulants and their
  combinatorics, see \cite{L}.
\end{remark}
Equation \eqref{eq:cumulants} expresses the joint moment $\E(X_1\cdots
X_n)$ in terms of products of cumulants.  For the first few values of
$n$, this equation is given by:
\[\E(X_1)=\kappa_1(X_1)\]
\[\E(X_1X_2)=\kappa_2(X_1\otimes X_2)+\kappa_1(X_1)\kappa_1(X_2)\]
\begin{multline*}\E(X_1X_2X_3)=
  \kappa_3(X_1\otimes X_2\otimes X_3)+\kappa_1(X_1)\kappa_2(X_2\otimes
  X_3)\\+\kappa_2(X_1\otimes
  X_2)\kappa_1(X_3)+\kappa_1(X_1)\kappa_1(X_2)\kappa_1(X_3)
\end{multline*}
One easily solves these equations for the cumulants expressed in terms
of products of joint moments:
\[\kappa_1(X_1)=\E(X_1)\]
\[\kappa_2(X_1\otimes X_2)=\E(X_1X_2)-\E(X_1)\E(X_2)\]
\begin{multline*}
  \kappa_3(X_1\otimes X_2\otimes
  X_3)=\E(X_1X_2X_3)-\E(X_1X_2)\E(X_3)\\-\E(X_1)\E(X_2X_3)+\E(X_1)\E(X_2)\E(X_3).
\end{multline*}
and in general finds
\begin{equation}\label{eq:moments}
  \kappa_n=\sum_{k=1}^n
  (-1)^{k-1} \alpha'_k\circ \left(\sum_P \bigotimes_{j=1}^k \E\circ \alpha_{n_j}\right).
\end{equation}
The sum in Equation \eqref{eq:moments} above is over the same set: $P$
ranges over all ordered partitions $P=(n_1,\ldots, n_k)$ with
$\displaystyle\sum_{j=1}^k n_j = n$.
\subsection{Main Proposition}
\begin{hypotheses*}
  Let $(V,\E,a)$ be a probability space.  Consider both $V$ and the
  complex numbers $\C$ as graded vector spaces concentrated in degree
  zero. Then $(V,0)$ and $(\C,0)$ are $A_\infty$ algebras and the map
  $\E:V\to \C$ defines an $A_\infty$ morphism between these two
  $A_\infty$ algebras.  Denote by $E:TV\to T\C$ the lift of $\E$ as a
  coalgebra map.  Furthermore, following the notation of
  Section~{\ref{sec:gauge}}, let $\mathbf{a}=\exp(A)$ where $A:TV\to TV$ is
  the lift of $a:V\otimes V \to V$ as a coderivation, and similarly let
  $\mathbf{a'}=\exp(A')$, where $A':T\C \to T\C$ is the lift of $a':\C
  \otimes \C \to \C$ to a
  coderivation.  Then,
  $E^{\mathbf{a},\mathbf{a}'}:TV\to T\C$ is
  an $A_\infty$ morphism between $(V,0^{\mathbf{a}})=(V,0)$ and
  $(\C,0^{\mathbf{a}'})=(\C,0)$.  Let $e_n:V^{\otimes n}\to \C$ denote
  the components of the $A_\infty$ morphism
  $E^{\mathbf{a},\mathbf{a}'}$.
\end{hypotheses*}
\begin{mainproposition*} \label{prop1} $\kappa_n=e_n$ for all $n$.
\end{mainproposition*}
\begin{proof}
  Both the collection $\{\kappa_n:V^{\otimes n}\to \C,n\ge 1\}$ and
  the single map $\E:V\to \C$ can be extended as coalgebra maps $TV\to
  T\C$. Let $K$ and $E$, respectively, denote these extensions. The
  statement of the proposition is that these two coalgebra maps are
  related by means of the coalgebra isomorphisms $\ba$ and $\ba'$ as
  in the following diagram.
  \[
  \begin{CD}
    TV @> K >> T\C\\
    @V\mathbf{a}VV @VV\mathbf{a}'V\\
    TV @> E >> T\C
  \end{CD}
  \]

  It suffices to check that the components
  \begin{gather*}
    \begin{CD}
      TV @> \mathbf{a} >> TV @> E >> T\C @> >> \C
    \end{CD}\\
    \begin{CD}
      TV @> K >> T\C @> \mathbf{a}' >> T\C @> >> \C
    \end{CD}
  \end{gather*}
  coincide when evaluated on a vector $X_1 \otimes \cdots\otimes X_n
  \in V^{\otimes n}$.

  The map $TV \stackrel{E}{\to} T\C\to \C$ is zero except for the one
  to one component $\E:V\to \C$.  Recall that $A$ denotes the lift of
  $a$ to a coderivation $TV\to TV$ and $\mathbf{a}=\exp(A)$. So the
  only components of $\mathbf{a}$ that contribute to the composition
  in question are $\frac{1}{(n-1)!}  A^{n-1}:V^{\otimes n}\to V$.

  The nonzero components of $A$ are of the form $a_n:V^{\otimes n}\to
  V^{\otimes n-1}$ and are given by
  \[
  a_n(X_1\otimes\cdots \otimes X_n)=\sum_{j=1}^{n-1} X_1\otimes\cdots
  \otimes X_jX_{j+1}\otimes\cdots\otimes X_n.\] Then the only nonzero
  component of $\frac{1}{(n-1)!}A^{n-1}$ is the composition
  $\frac{1}{(n-1)!}a_2\circ\cdots\circ a_n$. The expression for $a_n$
  is a sum of $(n-1)$ terms, so this composition has $(n-1)!$ terms. Applied
  to $X_1\otimes\cdots \otimes X_n$, this composition then yields
  $X_1\cdots X_n$. So the composition along the left and bottom takes
  $X_1\otimes \cdots \otimes X_n$ to \[\E(X_1\cdots X_n).\]

  The map $K:TV \to T\C$ evaluated on $(X_1\otimes\cdots\otimes X_n)$
  breaks into the following sum.
  \begin{equation}
    K(X_1\otimes\cdots\otimes X_n) = \left(\sum_{k=1}^n\sum_P
      \bigotimes_{j=1}^k \kappa_{n_j}\right)(X_1\otimes\cdots\otimes X_n)
  \end{equation}
  where $P$ ranges over all ordered partitions $P=(n_1,\ldots, n_k)$
  with $\displaystyle\sum_{j=1}^k n_j = n$.

  As in the calculation above for $\mathbf{a}$, the nonzero component
  of $\mathbf{a}'$ mapping $\C^{\otimes k}\to \C$ is
  $\frac{1}{(k-1)!}(A') ^{k-1}$ which maps $z_1\otimes \cdots \otimes
  z_n \mapsto z_1\cdots z_n$.  Hence, the composition along the top
  and right yields
  \begin{equation}
    \sum_{k=1}^n
    \alpha' \circ \left(\sum_P\bigotimes \kappa_{n_j}\right)
  \end{equation}
  where the sum is over the same partitions as above.  The coincidence
  of the two maps now follows from Definition \ref{def:cumulants} of
  the cumulants.
\end{proof}

\begin{remark}
  In classical probability theory, random variables are measurable
  $\C$-valued functions on a measure space and the expectation value
  of a random variable is defined by integration.  The product of
  measurable functions is measurable and defines the product of random
  variables.  In this situation, the product is commutative and
  associative.  One can define a \emph{classical probability space} as
  a probability space $(V,\E,a)$ for which $a:V\times V \to V$ is
  commutative.  The entire discussion in Sections \ref{sec:shadow} and
  \ref{sec:cumulants} of this paper can be symmetrized for a classical probability
  space.  The requisite modifications and results are contained in
  Section 2 of \cite{HPT2}.
\end{remark}

\section{Homotopy Probability Theory}
The starting point of homotopy probability theory is to replace the
space $V$ of random variables with a chain complex $C=(V,d)$ of random
variables.

\begin{definition}The data of a \emph{homotopy probability space}
  consists of a chain complex $C=(V,d)$, a chain map $\E:C \to \C$, and
  a degree zero associative product $a:V^{\otimes 2} \to V$.
\end{definition}

The expectation $\E$ and the differential $d$ are not assumed to
satisfy any properties with respect to $a:V^{\otimes 2}\to V$.

The coincidence of the cumulants for a probability space and an
$A_\infty$ morphism on the expectation provides the guide for how to
proceed for a homotopy probability space.  Cumulants for a homotopy
probability space are defined as the $A_\infty$ morphism $E^{\ba,
  \ba'}$ on the chain map $\E:(V,d)\to (\C,0)$ associated to the
product $a:V^{\otimes 2}\to V$ and the product $a':\C^{\otimes 2}\to
\C$ of complex numbers.  These cumulants are an $A_\infty$ morphism
between the $A_\infty$ structure $(V,D^\ba)$ and the (zero) $A_\infty$
structure $(\C,D^{\ba'})$.

These ideas are elaborated in Sections 2 and 3 of \cite{HPT2} and
illustrated with an example in Section 4 of \cite{HPT2}.

\nocite{Park2011}

\end{document}